\documentclass[11pt,draft,amsmath,amscd,amsbsy,amssymb]{amsart}
 \font\fivrm=cmr5
\relax

\chardef\bslash=`\\ % p.  424, TeXbook %\newcommand{\ntt}{\seriesm\shape n\tt}
\newcommand{\cs}[1]{{\protect\ntt\bslash#1}}
\newcommand{\opt}[1]{{\protect\ntt#1}} \newcommand{\env}[1]{{\protect\ntt#1}}
\makeatletter \def\verbatim{\interlinepenalty\@M \@verbatim
\leftskip\@totalleftmargin\advance\leftskip2pc
\frenchspacing\@vobeyspaces \@xverbatim} \makeatother \hfuzz1pc

\makeatletter \def\dgt@k{\dg@DX=-3 \dg@DY=2 \dg@SIZE=3}
\makeatother

\makeatletter \def\dgt@kk{\dg@DX=3 \dg@DY=-1 \dg@SIZE=3}% \makeatother

\theoremstyle{plain} \newtheorem{thm}{Theorem}[section]

\newtheorem{lemma}[thm]{Lemma}
\newtheorem{prop}[thm]{Proposition}

\theoremstyle{definition} 
 \newtheorem{que}[thm]{Question}

\newcommand{\R}{\mathcal R}
\newcommand{\RPlus}{\Real^{+}}
\newcommand{\norm}[1]{\left\Vert#1\right\Vert}
\newcommand{\abs}[1]{\left\vert#1\right\vert}
\newcommand{\set}[1]{\left\{#1\right\}}
\newcommand{\seq}[1]{\left<#1\right>}
\newcommand{\eps}{\varepsilon}
\newcommand{\To}{\longrightarrow}
\newcommand{\BX}{\mathbf{B}(X)}
\newcommand{\cw}{\operatorname{\mathrm{cw}}}
\newcommand{\pcw}{\operatorname{\mathrm{pcw}}}
\newcommand{\crw}{\operatorname{\mathrm{crw}}}
\newcommand{\cc}{\operatorname{\mathrm{cc}(\mathbb R^2)}}
\newcommand{\M}{\mathcal{M}}
\newcommand{\N}{\mathcal{N}}
\newcommand{\Lom}{\mathcal{L}}
\newcommand{\Comp}{\mathcal{K}}
\newcommand{\pr}{\mathrm{pr}}

\newcommand{\CCI}{\mathcal{I}}

\newcommand{\CCA}{\mathcal{A}}

\usepackage[all]{xy}
\begin{document}

\begin{abstract}
\end{abstract}

\title[]
{Gromov-Hausdorff ultrametric}

\author{I. Zarichnyi}
\address{Department of Mechanics and Mathematics, Lviv National University}
\email{izar@litech.lviv.ua}
%\thanks{}

%\date{}
%\dedicatory{}
%\commby{}

%%% ----------------------------------------------------------------------

\begin{abstract} We show that there exists a natural counterpart
of  the Gromov-Hausdorff metric in the class of ultrametric
spaces. It is proved, in particular, that the space of all
ultrametric spaces whose metric take values in a fixed countable
set is homeomorphic to the space of irrationals.
\end{abstract}

%%% ----------------------------------------------------------------------
\maketitle
%%% ----------------------------------------------------------------------

\section{Introduction}

Recall that the {\it Hausdorff distance} between two nonempty
closed bounded subsets, $A$ and $B$, of a metric space is
evaluated by the formula
$$d_H(A,B)=\inf\{\varepsilon>0\mid A\subset O_\varepsilon(B),\
B\subset O_\varepsilon(A)\}.$$

Given two compact metric spaces, $(X, d_X)$ and $(Y,d_Y)$,
 the Gromov-Hausdorff distance between them is defined by the formula
$$\varrho_{GH}(X,Y)=\inf\{d_H(i(X),j(Y))\mid i\colon X\to Z,\ j\colon Y\to Z
\text{ are isometric embeddings }\}.$$

Recall that a metric $d$ on $X$ is called an {\it ultrametric} if
it satisfies the following strong triangle inequality:

$$d(x,y)\le \max\{d(x,z),d(z,y)\},\ x,y,z\in X.$$

We are going to define a version of the Gromov-Hausdorff distance
for ultrametric spaces.

Given two compact ultrametric spaces, $(X, d_X)$ and $(Y,d_Y)$, we
define
\begin{align*}\varrho_{GHu}(X,Y)=&\inf\{d_H(i(X),j(Y))\mid i\colon X\to Z,\
j\colon Y\to Z\\ &\text{are isometric embeddings, where }Z\text{
is an ultrametric space}\}.\end{align*}

One can easily see that $\inf$ is well-defined as for every two
ultrametric spaces there exists an ultrametric space in which they
can be isometrically embedded.

\begin{lemma}\label{l:1}  Let $X_1,X_2$ be ultrametric spaces, $X_1\cap X_2=A$ and the
restrictions of ultrametrics in $X_1$ and $X_2$ onto $A$ coincide.
Then the formula
$$d(x_1,x_2)=\inf\{\max\{d_1(x_1,a),d_2(a,x_2)\}\mid a\in A\},$$
together with the initial ultrametrics on $X_1$ and $X_2$,
determines an ultrametric on $X_1\cup X_2$.
\end{lemma}
\begin{proof} We are going to prove the strong triangle inequality. Let $x,y,z\in X=X_1\cup X_2$. Without loss of
generality, one may assume that $x,y\in X_1\setminus X_2$, $z\in
X_2\setminus X_1$. There exist $a,b\in A$ such that $$d(x,z)=
\max\{d_1(x,a),d_2(a,z)\},\ d(y,z)= \max\{d_1(y,b),d_2(b,z)\}.$$

For the sake of brevity, we introduce the following notations:
 \begin{align*}
 \alpha=d_1(a,x),&\ \beta=d_1(b,y),\\ \gamma=d_2(a,z),&\
 \delta=d_2(b,z),\\ \eta=d_1(x,y),&\ \epsilon=d(x,z),\
 \zeta=d(z,y).\end{align*}

 The rest of the proof consists in analyzing all possible cases.

First, we are going to show that $d(x,y)\le
\max\{d(x,z),d(z,y)\}$, i.e. $\eta\le\max\{\epsilon,\zeta\}$.

1) $\epsilon=\alpha$, $\zeta=\beta$.  Suppose, on the contrary,
that $\eta>\max\{\epsilon,\zeta\}$, then $d_1(x,b)=\eta$,
$d_1(a,b)=\eta$ and, therefore,
$\eta\le\max\{\gamma,\delta\}\le\max\{\alpha,\beta\}<\eta$, a
contradiction.

2)  $\epsilon=\alpha$, $\zeta=\delta$.  Suppose that
$\eta>\max\{\epsilon,\delta\}$, then $d_1(a,y)=\eta$ and, since
$\eta>\delta\ge\beta$, we see that $d_1(a,b)=\eta$. Thus
$\alpha\ge\gamma=\eta$, a contradiction.

3)  $\epsilon=\alpha$, $\zeta=\delta$. Suppose that
$\eta>\max\{\epsilon,\zeta\}\ge\max\{\gamma,\delta\}$. Since
$\alpha\le\gamma<\eta$, we have $d_1(a,y)=\eta$. In turn, since
$\beta\le\zeta<\eta$, we have $d_1(a,b)=d_2(a,b)=\eta$. But then
$$\eta=d_2(a,b)\le
\max\{\gamma,\delta\}\le\max\{\epsilon,\zeta\}<\eta$$ and we come
to a contradiction.

The case $\epsilon=\gamma$, $\zeta=\beta$ is treated similarly to
case 2).

Now we are going to show that $d(z,y)\le \max\{d(x,z),d(x,y)\}$,
i.e. $\zeta\le\max\{\epsilon,\eta\}$.

1)  $\epsilon=\alpha$, $\zeta=\beta$. Suppose, on the contrary,
that $\zeta>\max\{\epsilon,\eta\}$, then $\beta>\eta$ and
$d_1(x,b)=\beta=\zeta$. Since $\zeta>\alpha$, we see that
$d_1(a,b)=\beta=\zeta$. Since $\zeta>\alpha\ge\gamma$, we see that
$\delta=\zeta$. We have $d_1(a,y)\le\max\{\alpha,\eta\}<\zeta$.
Also $\gamma\le\alpha<\zeta$ and therefore we obtain a
contradiction $\zeta\le\max\{d_1(a,y),\gamma\}<\zeta$.

2)  $\epsilon=\alpha$, $\zeta=\delta$. Suppose that
$\zeta>\max\{\epsilon,\eta\}\ge\max\{\alpha,\eta\}$. We have
$\delta>\max\{\alpha,\eta\}\ge\gamma$ and, therefore,
$d_2(a,b)=d_1(a,b)=\delta$. Since $d_1(a,b)=\delta>\alpha$, we see
that $d_1(x,b)=\delta$.  We have
$$d_1(a,y)\le\max\{\eta,\alpha\}\le\max\{\eta,\epsilon\}<\delta.$$
Thus,
$\zeta\le\max\{d_1(a,y),\gamma\}\le\max\{d_1(a,y),\epsilon\}<\zeta$,
a contradiction.

3)  $\epsilon=\gamma$, $\zeta=\beta$. Suppose that
$\zeta>\max\{\epsilon,\eta\}=\max\{\gamma,\eta\}$. We have
$d_1(a,y)\le\max\{\alpha,\eta\}\le\max\{\epsilon,\eta\}<\zeta$.
Since $\gamma\le\epsilon<\zeta$, we obtain
$\zeta\le\max\{d_1(a,y),\gamma\}<\zeta$, a contradiction.

4)  $\epsilon=\gamma$, $\zeta=\delta$. Suppose that
$\zeta>\max\{\epsilon,\eta\}$. Then $d_2(a,b)=d_1(a,b)=\delta$.
Since
$d_1(a,y)\le\max\{\alpha,\eta\}\le\max\{\epsilon,\eta\}<\zeta$, we
obtain a contradiction
$\zeta\le\max\{d_1(a,y),\gamma\}\le\max\{d_1(a,y),\epsilon\}<\zeta$.

That $d(z,x)\le \max\{d(y,z),d(x,y)\}$ can be proven similarly.
\end{proof}
\begin{thm}
The function $\varrho_{GHu}$ is an ultrametric on the set of
isometry classes of ultrametric spaces.
\end{thm}

\begin{proof} The symmetry is obvious. Since $\varrho_{GH}\le
\varrho_{GHu}$, we see that $\varrho_{GHu}(A,B)>0$ for
nonisometric $A$ and $B$.

We are going to prove the strong triangle inequality. Let
$X_1,X_2,X_3$ be ultrametric spaces and let $\varepsilon>0$ be
given. There exist ultrametric spaces $Y$ and $Z$ and isometric
embeddings $i_k\colon X_k\to Y$, $k=1,2$ and $j_l\colon X_l\to Z$,
$l=2,3$, such that $$d_H(i_1(X_1),i_2(X_2))\le
\varrho_{GHu}(X_1,X_2)+\varepsilon,\ d_H(j_2(X_2),j_3(X_3))\le
\varrho_{GHu}(X_2,X_3)+\varepsilon.$$  Identify $i_2(X_2)$ with
$j_2(X_2)$ along the map $j_3i_2^{-1}$. We obtain the quotient
set, which we denote by $K$,  of the disjoint union $Y\sqcup Z$.
For the sake of notational simplicity, we naturally identify $Y$
and $Z$ with the subspaces of $K$. By Lemma \ref{l:1}, there
exists an ultrametric, $d$, on $K$ which extends initial
ultrametrics on $Y$ and $Z$. Since the Hausdorff metric on the
space of nonempty compact subsets of an ultrametric space is an
ultrametric, we see that, in $K$,
$$d_H(i_1(X_1),j_3(X_3)\le\max
\{d_H(i_1(X_1),i_2(X_2)),d_H(i_2(X_3),j_3(X_3))\}$$ and therefore
\begin{align*}\varrho_{GHu}(X_1,X_3)\le& d_H(i_1(X_1),j_3(X_3)\le\max
\{d_H(i_1(X_1),i_2(X_2)),d_H(i_2(X_3),j_3(X_3))\}\\ \le&
\max\{\varrho_{GHu}(X_1,X_2)+\varepsilon,\
\varrho_{GHu}(X_2,X_3)+\varepsilon\}.
\end{align*}
Tending $\varepsilon$ to 0, we are done.
\end{proof}
\section{Ultrametric Gromov-Hausdorff space}

By $U$ we denote the  Gromov-Hausdorff space, i.e. the space of
all isometry classes of compact ultrametric spaces endowed with
the Gromov-Hausdorff ultrametric. For the sake of simplicity, we
prefer to work with representatives of the isometry classes rather
than with the classes themselves.

Denote by $\exp X$ the set of all nonempty compact subsets in $X$
endowed with the Hausdorff metric. It is well-known (see, e.g.,
\cite{ff}) that $\exp X$ is complete if so is $X$.

\begin{prop} The space $U$ is complete.
\end{prop}
\begin{proof} Let $(X_i)_{i=1}^\infty$ be a Cauchy sequence in
$U$. Without loss of generality, one may assume that $X_i$ and
$X_{i+1}$ lie in the same ultrametric space, $Y_i$. Let $Y=\sqcup
\{Y_i\mid i\in\mathbb N\}$. Similarly as in the proof of Lemma
\ref{l:1}, we subsequently glue $Y_2$ to $Y_1$ along $X_1$, then
glue the resulting space to $Y_3$ along $X_2$ etc. We obtain the
expanding sequence of  ultrametric spaces $Y_1$,
$Y_1\sqcup_{X_2}Y_2$, $Y_1\sqcup_{X_2}Y_2\sqcup_{X_3}Y_3$, \dots.
Let $Y$ denote the union of this sequence. Obviously, $Y$ is an
ultrametric space and therefore so is its completion, which we
denote by $\tilde Y$. The spaces $X_i$ are naturally embedded into
$\tilde Y$ and the sequence $(X_i)$ is a Cauchy sequence in
$\tilde Y$. Since the space $\exp \tilde Y$ is complete, there
exists the limit of the sequence $(X_i)$ in this space, which we
denote by $X$. It is evident that $X$ is also the limit of the
sequence $(X_i)$ in the space $U$.

\end{proof}
\begin{prop} The space $U$ is not separable.
\end{prop}
\begin{proof} For any $c\in [1/2,1]$, denote by $X_c$ the
two-point metric space with the nonzero distance equal to $c$. We
are going to prove that $\varrho_{GHu}(X_{c_1},X_{c_2})\ge1/4$
whenever $c_1\neq c_2$. Indeed, otherwise one can embed $X_{c_1}$
and $X_{c_2}$ in some ultrametric space so that the Hausdorff
distance between the images is $<1/2$. Without loss of generality
one may assume that there is an ultrametric, $d$,  on the union
$X_{c_1}\cup X_{c_2}$ extending the initial ultrametrics on
$X_{c_1}=\{a_1,a_2\}$ and  $X_{c_2}=\{b_1,b_2\}$ and
$d(a_1,b_1)<1/4$, $d(a_1,b_1)<1/4$. It follows from the strong
triangle inequality that
$c_1=d(a_1,a_2)=d(a_1,b_2)=d(b_1,b_2)=c_2$ and we obtain a
contradiction.
\end{proof}

Given a subset $K\subset\mathbb R_+$ with $0\in K$, we denote by
$U(K)$ the set of all ultrametric spaces $(X,d)$ with $d(X\times
X)\subset K$.

\begin{lemma} The space $U(K)$ is a closed subspace of $U$, for
any $K\subset\mathbb R_+$ with $0\in K$.
\end{lemma}
\begin{proof} Let $(X_i)_{i=1}^\infty$ be a sequence in
$U(K)$ converging to $X\in U$. Assume, on the contrary, that
$X\notin U(K)$, then there exist $a,b\in X$ such that
$d(a,b)\notin K$. There exists $i$ such that $\varrho_{GH}(X,X_i)<
\frac12 d(a,b)$. Without loss of generality, one may assume that
$X,X_i$ are subsets of an ultrametric space $Z$ with
$d_H(X,X_i)<\frac12 d(a,b)$. There exist $a'b'\in X_i$ such that
$d(a,a')<\frac12 d(a,b)$, $d(b,b')<\frac12 d(a,b)$. It follows
from the triangle $a,b,a'$ that $d(a',b)=d(a,b)$. Similarly, it
follows from the triangle $a',b,b'$ that
$d(a',b')=d(a',b)=d(a,b)$. We obtain a contradiction with the fact
that $X_i\in U(K)$.
\end{proof}

\begin{thm} Let $K$ be a countable subset of $\mathbb R_+$ with
$0$ as its nonisolated point. Then the space $U(K)$ is
homeomorphic to the space of irrationals.
\end{thm}
\begin{proof} First of all note that the space $U(K)$ is
separable. To this end, we are going to demonstrate that the space
$U_f(K)=\{Y\in U(K)\mid |Y|<\infty\}$, which is easily  seen to be
countable,  is dense in $U(K)$.

Prove that  $U(K)$ is nowhere locally compact. Let $X\in U(K)$ and
$\varepsilon>0$. Consider a finite $\varepsilon$-net
$Y=\{x_1,\dots,x_k\}$ in $X$. Without loss of generality, we may
assume that $d(x_1,x_2)=\min\{d(x,y)\mid x,y\in X,\ x\neq y\}$.
There exists a positive $c\in K$ such that
$c<\min\{d(x_1,x_2),\varepsilon/2\}$. For every natural $n$,
define a metric space $Y_n$ as follows. Let
$Y_n=Y\sqcup\{1,\dots,n\}$ and the metric $\varrho$ on $Y$ is
defined by the conditions $\varrho|(Y\times Y)=d|(Y\times Y)$,
$\varrho(y,i)= d(y,x_1)$, for any $y\in Y$, $1\le i\le n$, and
$\varrho(i,j)=c$, for every $i,j\in\{1,\dots,n\}$, $i\neq j$.

An easy verification that $\varrho$ is an ultrametric on $Y_n$ is
left to the reader.

Next, we note that $\varrho(Y_m,Y_n)\le c$ for every $m, n$. In
addition,  from the pigeon hole principle it easily follows that
$\varrho(Y_m,Y_n)\ge c/2$, whenever $m\neq n$. Therefore, the set
$\{Y_i\mid i\in \mathbb N\}$ is a countable discrete subset of a
closed  $c$-neighborhood of $X$ in $U(K)$. This demonstrates that
the space $U(K)$ is nowhere locally compact.

Remark also that the space $U(K)$ being a closed subset of $U$ is
complete.

It follows from \cite{k} that the space $U(K)$ is homeomorphic to
the space of irrationals.

\end{proof}
\section{Open problems}

\begin{que} Describe the topology of the space $U$.
\end{que}

A generalization of ultrametric spaces is introduced by David and
Semmes \cite{DS}. A metric space $(X, d)$ is said to be {\it
uniformly disconnected} if there exists $c > 0$ such that
$\max\{d(x_i, x_{i-1})\mid i=1,\dots, N\}\le cd(x, y)$  for all
finite chains of points $x = x_0, x_1,\dots , x_N = y$. In
\cite{DS}  it is proved  that, for any metric space $(X, d)$, the
metric $d$ is bi-Lipschitz equivalent with an ultrametric on $X$
if and only if the space $(X, d)$ is uniformly disconnected. This
result allows to find a counterpart of  the notion of the
Gromov-Hausdorff metric in the class of uniformly disconnected
spaces.

\begin{que} Is the obtained space of compact uniformly disconnected
spaces separable?
\end{que}

It is proved in \cite{EPW} that the space of (rooted) compact real
trees is complete. Here it is assumed that the set of these trees
is endowed with the Gromov-Hausdorff metric. Like in the case of
ultrametric spaces, we obtain another metric if we restrict
ourselves with embeddings in trees. We conjecture that the analogy
between trees and ultrametric space (see, e.g., \cite{h}) can be
extended also to the case of the obtained hyperspaces).


\begin{thebibliography}{99}

\bibitem{G} M. Gromov, Metric Structures for Riemannian and Non-Riemannian Spaces.
(Translation by S.M. Bates of 1981 French edition).- Birkh\"auser:
Boston, 1999.
\bibitem{ff} J. Munkres, Topology (2nd edition).- Prentice Hall, 1999.

\bibitem{k}  K. Kuratowski, Topology. Vol. I, New edition, revised and augmented.
Translated from the French by J. Jaworowski.- New York: Academic
Press, 1966.


\bibitem{DS}  G. David, S. Semmes, Fractured Fractals and Broken Dreams: Self-Similar
Geometry through Metric and Measure. - Oxford Lecture Ser. Math.
Appl. 7, Oxford University Press, 1997.


\bibitem{EPW} S.N. Evans,  J.W. Pitman,   A. Winter, {\it Rayleigh processes,
real trees and root growth with re-grafting}. Probab. Th. Rel.
Fields, (2003) to appear.

\bibitem{h} B. Hughes, {\it Trees and ultrametric spaces: a categorical
equivalence}, Adv. Math. 189, N 1 (2004), 148--191.
\end{thebibliography}
\end{document}